\documentclass[10pt]{article}         

\usepackage{amsmath,amsfonts,amssymb}   
\usepackage{bm}
\usepackage[utf8]{inputenc}
\usepackage{listings}
\usepackage{cite}
\usepackage[title]{appendix}
\usepackage{graphicx}
\usepackage[labelfont=bf]{caption}
\usepackage{listings,xcolor}

\title{On The Spherical Clothoid}
\author{Alexandru Ionu\textcommabelow{t}}

\begin{document}

\maketitle

\begin{abstract}
We revisit a nonlinear spline primitive for 3-space first studied by Even Mehlum. It is the spherical clothoid, the spherical curve with geodesic curvature a linear function of arc length. We present its Cartesian coordinate functions using confluent hypergeometric functions (the Kummer functions) and its stereographic projection onto the complex plane. New Humbert series results are also presented along with generating function formulas related to the associated Meixner-Pollaczek polynomials.
\end{abstract}

\section{Notation}
 
The Pochhammer symbol:\\
\[ (q)_n=q (q+1) \dots (q+n-1)=\frac{\Gamma(q+n)}{\Gamma(q)} \]
Gauss' hypergeometric series \cite{gauss1813disquisitiones}: 
\[ _2F_1(a,b,c;z)= \sum_{m=0}^{\infty} \frac{(a)_m (b)_m}{(c)_m} \frac{z^m}{m!} \]
The confluent hypergeometric functions (the Kummer function)  \cite{Kummer+1837+228+242}:
\[ _1F_1(a,b;z)= \sum_{m=0}^{\infty} \frac{(a)_m}{(b)_m} \frac{z^m}{m!} \]
Humbert hypergeometric series in 2 variables \cite{zbMATH02606592, bateman1953higher}:
\[ \phi_1(a,b,c;x,y)= \sum_{m,n=0}^{\infty} \frac{(a)_{m+n} (b)_m}{(c)_{m+n}} \frac{x^m}{m!} \frac{y^n}{n!} \]
\[ \phi_2(b_1,b_2,c;x,y)= \sum_{m,n=0}^{\infty} \frac{(b_1)_m (b_2)_n}{(c)_{m+n}} \frac{x^m}{m!} \frac{y^n}{n!} \]
\[ \Xi_1(a_1,a_2,b,c;x,y)= \sum_{m,n=0}^{\infty} \frac{(a_1)_{m} (a_2)_{n} (b)_m}{(c)_{m+n}} \frac{x^m}{m!} \frac{y^n}{n!} \]
The parabolic cylinder functions, even and odd solutions to $\frac{d^2 y}{d z^2}-(\frac{1}{4}z^2+a)y=0$ (19.2.1, 19.2.2 in \cite{abramowitz+stegun}):
\[ y_1 (a;z) = e^{\frac{-z^2}{4}} \,_1F_1 \left( \frac{1}{2}a+\frac{1}{4},\frac{1}{2};\frac{z^2}{2} \right) \]
\[ y_2 (a;z) = z e^{\frac{-z^2}{4}} \,_1F_1 \left( \frac{1}{2}a+\frac{3}{4},\frac{3}{2};\frac{z^2}{2} \right)  \]
The associated Meixner-Pollaczek polynomials $Q_n^{\lambda}(x;\phi,c)$ \cite{poll, Luo2019GeneratingFO} defined by the recurrence relation:
\begin{align*} 
0 = &\phantom{+} (n+c+1) Q_{n+1}^{\lambda}(x;\phi,c) \\ 
    & -2[(n+\lambda+c)\cos \phi + x \sin \phi] Q_n^{\lambda}(x;\phi,c) \\ 
		& + (n+2 \lambda+c-1) Q_{n-1}^{\lambda}(x;\phi,c)
\end{align*}  
where $n=0,1,2...$ and $Q_{-1}^{\lambda}(x;\phi,c)=0$, $Q_{0}^{\lambda}(x;\phi,c)=1$
\section{Summary of prior work}

In mathematics, the word \emph{spline} has come to mean a function defined piecewise by polynomials with certain continuity constraints. We are interested in a more general concept i.e. nonlinear splines meaning the curve primitive, the building block of the spline, can be an analytic function. In English, The term \emph{spline} goes back centuries referring to the flexible strip of wood used by draftsmen to draw smooth curves. Lead weights called \emph{ducks} or \emph{whales} held the spline in place at certain points and the strip would take the shape that minimized the bending energy yielding a smooth curve. This shape arising in the physical world corresponds to the minimum energy curve (MEC) i.e. the curve that minimizes the $L^2$ norm of curvature (see \cite{Levien:EECS-2009-162} for a more comprehensive treatise).

We can also examine the MEC spline element in 3D space. First we recall that the intrinsic quantities needed to describe the shape of a space curve are its curvature $\kappa$ and torsion $\tau$ as functions of arc length $s$. Mehlum derived the relationship between curvature and torsion for this curve and a differential equation characterizing the curvature (1.4 and 1.3 in\cite{Meh94}):
\begin{equation}
\label{eqn: kappa tau}
	\kappa^2\tau = C
\end{equation}
\begin{equation}
\label{eqn: kappa diff EQ}
	\left [(\kappa^2)'\right ]^2 + \kappa^2\left [(\kappa^2 - 2D)^2 - 4\psi^2\right ] + 4C^2 = 0 
\end{equation}
$C$, $D$, and $\psi$ are constants. 

The calculus of variations problem is solved without the approximation that leads to the ubiquitous cubic spline. However, working with the MEC in practical applications can have its own drawbacks. Mehlum goes on to study a modified version of the differential equation for curvature. (1.5 and 1.6 in \cite{Meh94}) Assuming $\kappa^2 \ll 2D$ yields:
\begin{equation}
\label{eqn: spherical clothoid diff EQ}
	\left [(\kappa^2)' \right ]^2 - 4 \alpha^2 \kappa^2 + 4C^2 = 0 
\end{equation}
$\alpha=\psi^2 - D^2$ is a constant. We have a simple solution:
\begin{equation}
\label{eqn: spherical clothoid kappa squared}
	\kappa^2 = \alpha^2s^2 + \frac{C^2}{\alpha^2} 
\end{equation}

When $C=0$ the resulting planar curve is the clothoid i.e. curvature is a linear function of arc length. The clothoid has since been heavily investigated in the theory of nonlinear planar splines \cite{Stoer82,Coope93,Levien:EECS-2009-162}. For all other $C$, Mehlum discovered that this curve lies on a sphere and that $R=\frac{\alpha}{C}$ (2.16 in \cite{Meh94}). Mehlum and Wimp bring up a known necessary and sufficient condition for a space curve to be spherical (2.8 in \cite{mehlum_wimp_1985}):
\begin{equation}
\label{eqn: spherical kappa tau}
	\frac{\tau}{\kappa} - \Big ( \frac{\kappa'}{\tau\kappa^2} \Big )' = 0 
\end{equation}
Mehlum and Resch \cite{Meh94} go on to describe this curve with a geometric/kinematic construction by rolling a sphere without slipping or twisting on a planar clothoid. They thereby reveal that this is the spherical curve with geodesic curvature a linear function of arc length (see\cite{1534-0392_2014_1_435}, Theorem 2.3):
\begin{equation}
\label{eqn: kappa geee}
	\kappa_g = \alpha s 
\end{equation}
It is fitting to call it the sherical clothoid, a term first coined by \"{U}lo Lumiste \cite{10.1007/BFb0083636} in work unrelated to the theory of nonlinear splines.

Mehlum manages to express the cartesian coordinate functions of the spherical clothoid using Humbert series $\phi_1$, $\phi_2$ and some more common functions.

\section{Novel results}

\subsection{Simplified hypergeometric Representation}

The main purpose of this work is to present a simpler hypergeometric representation of the coordinate functions. Without loss of generality, we will restrict our attention to the case $R=1$ since all other curves can be obtained with simple scaling. Let $\bm{r}(s)=(x(s),y(s),z(s))$ be the parametric equation of our space curve. We consider the differential equation and initial conditions first examined by Mehlum (1.13, 2.17, 2.18, 2.19 and 2.21 in \cite{Meh94}):
\begin{equation}
\label{eqn: mehlum's 4th order DE with initial conditions}
	\begin{cases}
		\bm{r}''''(s)+(\alpha^2 s^2 + 1)\bm{r}''(s)+3 \alpha^2 s \bm{r}'(s)=\bm{0}\\
		\bm{r}(0)=(0,0,0)\\
		\bm{r}'(0)=(1,0,0)\\
		\bm{r}''(0)=(0,1,0)\\
		\bm{r}'''(0)=(-1,0,\alpha)\\
  \end{cases}
\end{equation}
The following solution holds:
\begin{equation}
\label{eqn: x,y,z solution}
	\begin{cases}
		x(s)=\Re \left( s \,_1F_1\left( -\frac{\mathrm{i}}{8 \, \alpha}, \frac{1}{2}; -\frac{\mathrm{i} \, \alpha s^{2}}{2} \right) \,_1F_1\left(\frac{\mathrm{i}}{8 \, \alpha} + \frac{1}{2}, \frac{3}{2} ; \frac{\mathrm{i} \, \alpha s^{2}}{2} \right) \right) \\
		y(s)=-2 {\left| \,_1F_1\left( \frac{\mathrm{i}}{8 \, \alpha} , \frac{1}{2}  ; \frac{\mathrm{i} \, \alpha s^{2}}{2} \right) \right|}^{2} + 2 \\
		z(s)=\Im \left( s \,_1F_1\left( -\frac{\mathrm{i}}{8 \, \alpha}, \frac{1}{2}; -\frac{\mathrm{i} \, \alpha s^{2}}{2} \right) \,_1F_1\left(\frac{\mathrm{i}}{8 \, \alpha} + \frac{1}{2}, \frac{3}{2} ; \frac{\mathrm{i} \, \alpha s^{2}}{2} \right) \right)
  \end{cases}
\end{equation}
$x$, $y$ and $z$ can also be expressed equivalently using conjugation instead of using real and imaginary parts and the absolute value. Consult Appendix \ref{appendix:a} for a computer-assisted proof using Mathematica. 

\begin{figure}[h]
\centering
\includegraphics[width=0.5\textwidth]{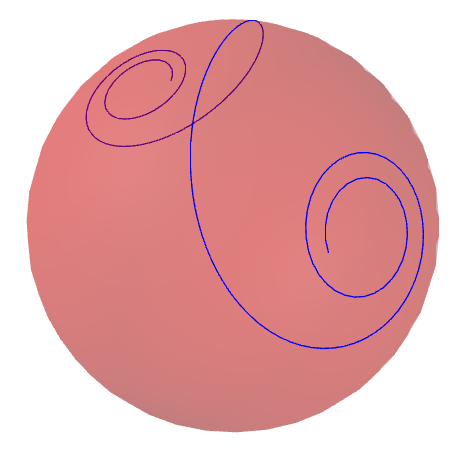}
\caption{Plot of the spherical clothoid for $\alpha=1$ and $s \in [-5,5]$}
\end{figure}

\subsection{Derivation}
We begin by investigating a loose end from the kinematic construction. The following system of equations describes the rolling without slipping and twisting of a unit sphere on a planar clothoid with scale factor $\alpha$ (see 3 in \cite{KholodenkoSilagadze} and 29a in \cite{RojoBloch}) : 
\begin{equation}
\label{eqn: 3d kinematics} 
	 \frac{d}{ds} \begin{pmatrix} \tilde{x} \\ \tilde{y} \\ \tilde{z} \end{pmatrix} = \begin{pmatrix} 0 & 0 & -\sin \left( \frac{\alpha s^2}{2} \right) \\ 0 & 0 & \cos \left( \frac{\alpha s^2}{2} \right) \\ \sin \left( \frac{\alpha s^2}{2} \right) & -\cos \left( \frac{\alpha s^2}{2} \right) & 0 \end{pmatrix} \begin{pmatrix} \tilde{x} \\ \tilde{y} \\ \tilde{z} \end{pmatrix}
\end{equation}
We introduce two new complex variables (see 10 in \cite{KholodenkoSilagadze}, 8 in \cite{RojoBloch}):
\begin{equation}
\label{eqn: Hopf map}
	 \tilde{x}=ab^*+ba^*,\tilde{y}= \mathrm{i} \left( ab^*-ba^* \right),\tilde{z}=aa^*-bb^*.
\end{equation} 
This tool was first presented by Feynman, Vernon and Hellwarth to geometrically represent the Schr\"{o}dinger equation of a two-level quantum system \cite{FeynmanVernonHellwarth}. $a$ and $b$ must satsfy the following system in order for $\tilde{x}$,$\tilde{y}$ and $\tilde{z}$ to satisfy (\ref{eqn: 3d kinematics}):
\begin{equation}
\label{eqn: complex system}
	 \mathrm{i} \frac{d}{ds} \begin{pmatrix} a \\ b \end{pmatrix} = -\frac{1}{2} \begin{pmatrix} 0 & e^{-\mathrm{i} \alpha s^2/2} \\ e^{\mathrm{i}\alpha s^2/2} & 0 \end{pmatrix} \begin{pmatrix} a \\ b \end{pmatrix}
\end{equation}
We are now left with a system that can be solved directly in a computer algebra system. Consult Appendix \ref{appendix:b} for a solution using Mathematica.  This system arises in the Landau-Zener problem (4, \cite{Zener}). In fact, Zener first solved this system analytically using parabolic cylinder functions in this quantum physics context. 
\\
We then obtain the solution to Mehlum's equation using (\ref{eqn: Hopf map}) again and a rigid motion to satisfy the initial conditions.

\subsection{Stereographic projection}

\begin{figure}[h]
\centering
\includegraphics[width=0.5\textwidth]{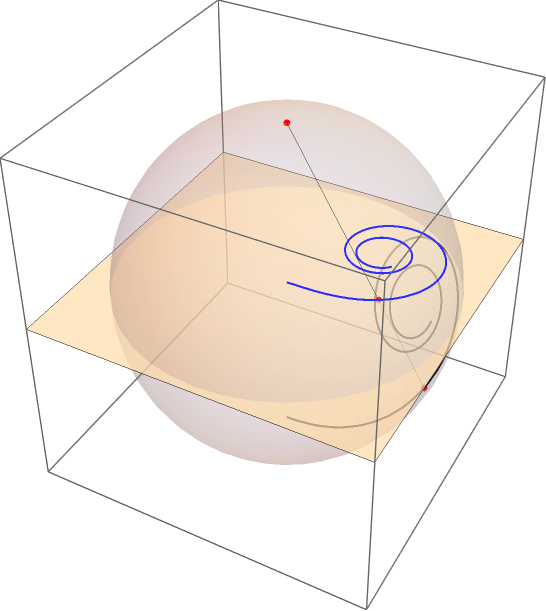}
\caption{Plot of the spherical clothoid and its stereographic projection for $\alpha=1$ and $s \in [0,5]$}
\end{figure}

Stereographic projection is a mapping that projects the sphere onto the plane. Circles on the sphere are mapped to circles on the plane (as long as they do not pass through the point of projection) and loxodromes are mapped to logarithmic spirals. We shift the sphere to center it on the origin and consider the stereographic projection of our curve onto the complex plane:
\begin{equation}
\label{eqn: stereo}
	\begin{split}
		\zeta(s) &= X(s)+\mathrm{i} Y(s) \\
		&= \frac{x(s)+\mathrm{i} z(s)}{1-(y(s)-1)} \\
		&= \frac{s \,_1F_1\left(\frac{\mathrm{i}}{8 \, \alpha} + \frac{1}{2}, \frac{3}{2} ; \frac{\mathrm{i} \, \alpha s^{2}}{2} \right)}{2 \,_1F_1\left( \frac{\mathrm{i}}{8 \, \alpha} , \frac{1}{2}  ; \frac{\mathrm{i} \, \alpha s^{2}}{2} \right)}	
	\end{split}
\end{equation}
We also have an alternative representation using a quotient of odd and even parabolic cylinder functions: 
\begin{equation}
\label{eqn: stereo2}
	\zeta(s) = \frac{e^{-\mathrm{i}\pi/4} \; y_2 \left( \frac{\mathrm{i}}{4 \alpha} -\frac{1}{2} ; e^{\mathrm{i}\pi/4} \sqrt{\alpha} s \right)}{2 \sqrt{\alpha} \;\;\; \, y_1  \left( \frac{\mathrm{i}}{4 \alpha} -\frac{1}{2} ; e^{\mathrm{i}\pi/4} \sqrt{\alpha} s \right) }
\end{equation}

\subsection{Humbert series corollaries}
Let us revisit Mehlum and Wimp's work on the spherical clothoid. We begin with Mehlum's hypergeometric expression for $y$ restricted to the unit sphere (2.54 in \cite{Meh94}):
\begin{equation}
\label{eqn: Mehlum y}
	y(s)=1-\phi_2 \left(\frac{\mathrm{i}}{4 \, \alpha},-\frac{\mathrm{i}}{4 \, \alpha},\frac{1}{2}; \frac{\mathrm{i} \, \alpha s^{2}}{2} ,-\frac{\mathrm{i} \, \alpha s^{2}}{2} \right) 
\end{equation}
Combining this work with our new representation yields a novel reduction formula for a special form of $\phi_2$:
\begin{equation}
\label{eqn: phi2 reduction}
	\phi_2 \left( a,-a,\frac{1}{2};x,-x \right) = 2 \,_1F_1\left( \frac{a}{2}, \frac{1}{2} ; x \right) \,_1F_1\left( -\frac{a}{2}, \frac{1}{2} ; -x \right) -1 
\end{equation}
We can use some suggestions in Mehlum's work to obtain an expression for $x + \mathrm{i} z$ using $\phi_1$ special forms and less mysterious functions (see Appendix \ref{appendix:b} for the derivation):
\begin{equation}
\label{eqn: Mehlum x+iz}
	\begin{split}
		&u= \frac{\sqrt{\pi } \cos (\pi\frac{\mathrm{i}}{4 \, \alpha}) \Gamma \left(-\frac{\mathrm{i}}{4 \, \alpha}-\frac{1}{2}\right)}{\frac{\mathrm{i}}{4 \, \alpha} \Gamma \left(-\frac{\mathrm{i}}{8 \, \alpha}\right)^2} \\
		&v = -\frac{\sqrt{\pi } \cos (\pi\frac{\mathrm{i}}{4 \, \alpha}) \Gamma \left(\frac{\mathrm{i}}{4 \, \alpha}-\frac{1}{2}\right)}{2 \Gamma \left(\frac{\mathrm{i}}{8 \, \alpha}+\frac{1}{2}\right)^2} \\
		&\phi_1= \phi_1 \left(\frac{\mathrm{i}}{4 \, \alpha}+1,\frac{\mathrm{i}}{4 \, \alpha},\frac{\mathrm{i}}{4 \, \alpha}+\frac{3}{2};\frac{1}{2},\frac{\mathrm{i} \, \alpha s^{2}}{2} \right) \\
		&\phi_1^*= \phi_1 \left(-\frac{\mathrm{i}}{4 \, \alpha}+1,-\frac{\mathrm{i}}{4 \, \alpha},-\frac{\mathrm{i}}{4 \, \alpha}+\frac{3}{2};\frac{1}{2},-\frac{\mathrm{i} \, \alpha s^{2}}{2} \right) \\
		&x(s)+\mathrm{i} \, z(s) = s \left(u \, e^{-\frac{\mathrm{i} \, \alpha s^{2}}{2}} \phi_1 + v \, e^{\frac{\mathrm{i} \, \alpha s^{2}}{2}} \phi_1^*  \right)
	\end{split}
\end{equation} 
Once again combining this with our new representation yields a new identity:
\begin{equation}
\label{eqn: phi1 identity}
    \begin{split}
        &\phantom{-} \frac{\sqrt{\pi } \cos (\pi  a) \Gamma \left(-a-\frac{1}{2}\right)}{a \Gamma \left(-\frac{a}{2}\right)^2} e^{-x} \phi_1 \left(a+1,a,a+\frac{3}{2};\frac{1}{2},x \right) \\ 
        &-\frac{\sqrt{\pi } \cos (\pi  a) \Gamma \left(a-\frac{1}{2}\right)}{2 \Gamma \left(\frac{a+1}{2}\right)^2} e^{x} \phi_1 \left(-a+1,-a,-a+\frac{3}{2};\frac{1}{2},-x \right) \\[0.3cm]
        &\qquad \qquad \qquad= \,_1F_1\left( -\frac{a}{2}, \frac{1}{2}; -x \right) \,_1F_1\left(\frac{a}{2} + \frac{1}{2}, \frac{3}{2} ;x \right) 
    \end{split} 
\end{equation}

It is interesting to note that these reductions demystify some quadratic relationships for special functions studied by Mehlum and Wimp (section 5 in \cite{mehlum_wimp_1985}).  

\subsection{Generating function formulas related to the associated Meixner-Pollaczek polynomials}
Mehlum found expressions for the parameter functions of the spherical clothoid involving associated Meixner-Pollaczek polynomials (2.38, 2.39 and 2.40 in \cite{Meh94}): 
\begin{equation}
\label{eqn: gen x}
	x(s)= \sum_{j=0}^{\infty} Q_j^0(\frac{1}{4 \alpha};\frac{\pi}{2},\frac{1}{2}) \frac{\left( \frac{-\alpha }{2} \right)^j s^{2j+1}}{(1)_j}
\end{equation}
\begin{equation}
\label{eqn: gen y}
	y(s)= \frac{1}{2} \sum_{j=0}^{\infty} Q_j^0(\frac{1}{4 \alpha};\frac{\pi}{2},1) \frac{\left( \frac{-\alpha }{2} \right)^j s^{2j+2}}{\left( \frac{3}{2} \right)_j}
\end{equation}
\begin{equation}
\label{eqn: gen z}
	z(s)= \frac{\alpha}{6} \sum_{j=0}^{\infty} Q_j^0(\frac{1}{4 \alpha};\frac{\pi}{2},\frac{3}{2}) \frac{\left( \frac{-\alpha }{2} \right)^j s^{2j+3}}{(2)_j}
\end{equation}
Generating functions related to this family of orthogonal polynomials arise in recent research \cite{LuoRainaZhao, AhbliMouayn, LuoRaina}. Combining these expressions with our new parameter functions, we obtain the following generating function results: 
\begin{equation}
\label{eqn: gen result 1}
	\begin{split}
		\sum_{j=0}^{\infty} Q_j^0(x;\frac{\pi}{2},\frac{1}{2}) \frac{t^j}{j!} = &\frac{1}{2} \Bigg( \,_1F_1\left( -\frac{\mathrm{i}x}{2}, \frac{1}{2}; \mathrm{i} t \right) \,_1F_1\left(\frac{\mathrm{i}x}{2} + \frac{1}{2}, \frac{3}{2} ; -\mathrm{i} t \right) \\
		&+\,_1F_1\left( \frac{\mathrm{i}x}{2}, \frac{1}{2}; -\mathrm{i} t \right) \,_1F_1\left(-\frac{\mathrm{i}x}{2} + \frac{1}{2}, \frac{3}{2} ; \mathrm{i} t \right) \Bigg)
	\end{split}
\end{equation}

\begin{equation}
\label{eqn: gen result 2}
	\sum_{j=0}^{\infty} Q_j^0(x;\frac{\pi}{2},1) \frac{t^j}{\left( \frac{3}{2} \right)_j} =
	\frac{\,_1F_1\left(\frac{\mathrm{i}x}{2},\frac{1}{2};-\mathrm{i}t\right)\,_1F_1\left(-\frac{\mathrm{i}x}{2},\frac{1}{2};\mathrm{i}t\right)-1}{2xt}
\end{equation}

\begin{equation}
\label{eqn: gen result 3}
	\begin{split}
		\sum_{j=0}^{\infty} Q_j^0(x;\frac{\pi}{2},\frac{3}{2}) \frac{t^j}{(j+1)!} = &\frac{3\mathrm{i}}{2t} \Bigg( \,_1F_1\left( -\frac{\mathrm{i}x}{2}, \frac{1}{2}; \mathrm{i} t \right) \,_1F_1\left(\frac{\mathrm{i}x}{2} + \frac{1}{2}, \frac{3}{2} ; -\mathrm{i} t \right) \\
		&-\,_1F_1\left( \frac{\mathrm{i}x}{2}, \frac{1}{2}; -\mathrm{i} t \right) \,_1F_1\left(-\frac{\mathrm{i}x}{2} + \frac{1}{2}, \frac{3}{2} ; \mathrm{i} t \right) \Bigg)
	\end{split}
\end{equation}

\section{Conclusion}
We have found simple expressions for the Cartesian coordinates of the spherical clothoid and its projection with some special function results as a bonus. It is remarkable that tools from quantum physics elucidate a classical problem. Suddenly this curve studied by Mehlum does not seem so exotic, heightening its potential in the world of computer aided geometric design. The computation of the special functions presented in this work presents an avenue for future research. 

\section{Acknowledgments}

I would like to thank Dr. Zurab Silagadze for discussions that sparked the idea behind this hypergeometric reduction tied to the spherical clothoid. I am grateful to Dr. Khalid Ahbli,
Dr. James Hateley, Max K\"{o}lbl and Dr. Robert Lewis for their attention and support during my research.

\begin{appendices}
\section{ODE solution verification}
\label{appendix:a}
  \begin{lstlisting}[language=Mathematica,mathescape,numbers=none,title={Mathematica Session}]
    In[1]:= x[s_]:=1/2 s Hypergeometric1F1[-(I/(8 $\alpha$)),1/2,-((I $\alpha$ s^2)/2)] Hypergeometric1F1[1/2+I/(8 $\alpha$),3/2,(I $\alpha$ s^2)/2]+1/2 s Hypergeometric1F1[I/(8 $\alpha$),1/2,(I $\alpha$ s^2)/2] Hypergeometric1F1[1/2-I/(8 $\alpha$),3/2,-((I $\alpha$ s^2)/2)]
		In[2]:= y[s_]:=2-2 Hypergeometric1F1[-(I/(8 $\alpha$)),1/2,-((I $\alpha$ s^2)/2)] Hypergeometric1F1[I/(8 $\alpha$),1/2,(I $\alpha$ s^2)/2]
		In[3]:= z[s_]:=1/(2 I) s Hypergeometric1F1[-(I/(8 $\alpha$)),1/2,-((I $\alpha$ s^2)/2)] Hypergeometric1F1[1/2+I/(8 $\alpha$),3/2,(I $\alpha$ s^2)/2]-1/(2 I) s Hypergeometric1F1[I/(8 $\alpha$),1/2,(I $\alpha$ s^2)/2] Hypergeometric1F1[1/2-I/(8 $\alpha$),3/2,-((I $\alpha$ s^2)/2)]
		In[4]:= r[s_]:={x[s],y[s],z[s]}
		In[5]:= r[0]
		Out[5]= {0,0,0}
		In[6]:= r'[0]
		Out[6]= {1,0,0}
		In[7]:= r''[0]
		Out[7]= {0,1,0}
		In[8]:= FullSimplify[r'''[0]]
		Out[8]= {-1,0,$\alpha$}
		In[9]:= FullSimplify[r''''[s]+($\alpha$^2 s^2 + 1) r''[s] + 3 $\alpha$^2 s r'[s]]
		Out[9]= {0,0,0}
  \end{lstlisting}

\pagebreak 
\section{Complex differential equation system}
\label{appendix:b}

Here we take an experimental approach and consider the following initial conditions: $a(0)=1, b(0)=0$ (implying $\tilde{x}(0)=0,\tilde{y} (0)= , \tilde{z} (0)=1$)

\begin{lstlisting}[language=Mathematica,mathescape,numbers=none,title={Mathematica Session}]
In[1]:=DSolve[{a'[s] == 1/2 I b[s] E^(1/2 (-I) $\alpha$ s^2), b'[s] == 1/2 I a[s] E^(1/2 I $\alpha$ s^2), a[0] == 1, b[0] == 0}, {a, b}, s]
Out[1]={{a -> Function[{s}, E^(-(1/2) I s^2 $\alpha$) Hypergeometric1F1[-((-I - 4 $\alpha$)/(8 $\alpha$)), 1/2, 1/2 I s^2 $\alpha$]], b -> Function[{s}, 1/2 s (I Hypergeometric1F1[1 - (-I - 4 $\alpha$)/(8 $\alpha$), 3/2, 1/2 I s^2 $\alpha$] + 4 $\alpha$ Hypergeometric1F1[1 - (-I - 4 $\alpha$)/(8 $\alpha$), 3/2, 1/2 I s^2 $\alpha$] - 4 $\alpha$ Hypergeometric1F1[-((-I - 4 $\alpha$)/(8 $\alpha$)) 1/2, 1/2 I s^2 $\alpha$])]}}
In[2]:=FullSimplify[E^(-(1/2) I s^2 $\alpha$) Hypergeometric1F1[-((-I - 4 $\alpha$)/(8 $\alpha$)), 1/2, 1/2 I s^2 $\alpha$]]
Out[2]=Hypergeometric1F1[-(I/(8 $\alpha$)), 1/2, -(1/2) I s^2 $\alpha$]
In[3]:=FullSimplify[1/2 s (I Hypergeometric1F1[1 - (-I - 4 $\alpha$)/(8 $\alpha$), 3/2, 1/2 I s^2 $\alpha$] + 4 $\alpha$ Hypergeometric1F1[1 - (-I - 4 $\alpha$)/(8 $\alpha$), 3/2, 1/2 I s^2 $\alpha$] - 4 $\alpha$ Hypergeometric1F1[-((-I - 4 $\alpha$)/(8 $\alpha$)), 1/2, 1/2 I s^2 $\alpha$])]
Out[3]= 1/2 I s Hypergeometric1F1[1/2 + I/(8 $\alpha$), 3/2, 1/2 I s^2 $\alpha$]
\end{lstlisting}

\section{Mehlum's $x+\mathrm{i} \, z$}
\label{appendix:c}
Combining Mehlum's expression for $x$ and $z$ (2.65 and 2.67 in \cite{Meh94}) and a Humbert series transform (2.61 and 2.62 in \cite{Meh94}), it is possible to obtain expressions of this form:
\begin{equation}
\label{eqn: u v}
	\begin{split}
		x(s)+\mathrm{i} \, z(s) &= s \left(u \, e^{-\frac{\mathrm{i} \, \alpha s^{2}}{2}} \phi_1 + v \, e^{\frac{\mathrm{i} \, \alpha s^{2}}{2}} \phi_1^*  \right) \\
		&= s \left(u \, \Xi_1 + v \, \Xi_1^*  \right)
	\end{split}
\end{equation}
$u$, $v$ do not depend on $s$ and these are the special forms of $\Xi_1$ we are playing with:
\begin{equation}
\label{eqn: Xi1}
	\begin{split}
		\Xi_1 &= \Xi_1 \left( \frac{\mathrm{i}}{4 \, \alpha}, \frac{1}{2}, \frac{\mathrm{i}}{4 \, \alpha}+1, \frac{\mathrm{i}}{4 \, \alpha}+\frac{3}{2} ;\frac{1}{2},-\frac{\mathrm{i} \, \alpha s^{2}}{2} \right) \\
		&= \sum_{m=0}^{\infty} \frac{(\frac{1}{2})_m}{(\frac{\mathrm{i}}{4 \, \alpha}+\frac{3}{2})_m} {}_2F_1(\frac{\mathrm{i}}{4 \, \alpha},\frac{\mathrm{i}}{4 \, \alpha}+1,\frac{\mathrm{i}}{4 \, \alpha}+\frac{3}{2}+m;\frac{1}{2}) \frac{(-\frac{\mathrm{i} \, \alpha s^{2}}{2})^m}{m!} \\
		\Xi_1^* &= \Xi_1 \left( -\frac{\mathrm{i}}{4 \, \alpha}, \frac{1}{2}, -\frac{\mathrm{i}}{4 \, \alpha}+1, -\frac{\mathrm{i}}{4 \, \alpha}+\frac{3}{2} ;\frac{1}{2},\frac{\mathrm{i} \, \alpha s^{2}}{2} \right) \\
		&= \sum_{m=0}^{\infty} \frac{(\frac{1}{2})_m}{(-\frac{\mathrm{i}}{4 \, \alpha}+\frac{3}{2})_m} {}_2F_1(-\frac{\mathrm{i}}{4 \, \alpha},-\frac{\mathrm{i}}{4 \, \alpha}+1,-\frac{\mathrm{i}}{4 \, \alpha}+\frac{3}{2}+m;\frac{1}{2}) \frac{(\frac{\mathrm{i} \, \alpha s^{2}}{2})^m}{m!} \\
	\end{split}
\end{equation}
Let's examine the special forms of $_2F_1$ that arise. We first recall some identities derived by Mitra (2.68 and 2.69 in \cite{Meh94}):
\begin{equation}
\label{eqn: 2f1(a,a+1,a+1/2,1/2}
	_2F_1 \left( a,a+1,a+\frac{1}{2};\frac{1}{2} \right)= \sqrt{\pi} \Gamma \left( a + \frac{1}{2} \right)
	\left[ \frac{1}{\Gamma \left(\frac{a}{2}+\frac{1}{2} \right)^2} + \frac{1}{\Gamma \left(\frac{a}{2}+1 \right)\Gamma \left( \frac{a}{2} \right)}  \right]
\end{equation}
\begin{equation}
\label{eqn: 2f1(a,a+1,a+3/2,1/2}
	_2F_1 \left( a,a+1,a+\frac{3}{2};\frac{1}{2} \right)= 2 \sqrt{\pi} \Gamma \left( a + \frac{3}{2} \right)
	\left[ \frac{1}{\Gamma \left(\frac{a}{2}+\frac{1}{2} \right)^2} - \frac{1}{\Gamma \left(\frac{a}{2}+1 \right)\Gamma \left( \frac{a}{2} \right)}  \right]
\end{equation}
Using the contiguous relations of $_2F_1$ we can deduce another relation:
\begin{equation}
\label{eqn: 2f1(a,a+1,a+5/2,1/2}
	_2F_1 \left( a,a+1,a+\frac{5}{2};\frac{1}{2} \right)= \frac{4}{3} \sqrt{\pi} \Gamma \left( a + \frac{5}{2} \right)
	\left[ \frac{1-4a}{\Gamma \left(\frac{a}{2}+\frac{1}{2} \right)^2} + \frac{1+4a}{\Gamma \left(\frac{a}{2}+1 \right)\Gamma \left( \frac{a}{2} \right)}  \right]
\end{equation}
Examining the first few terms in the Maclaurin series of $s \left(u \, \Xi_1 + v \, \Xi_1^*  \right)$ and our initial conditions, we can now find expressions for $u$ and $v$ explicitly in terms of the gamma function and more common functions. There are no terms with even powers of $s$ in our Maclaurin series and that satisfies the first and third initial conditions. We examine the second and fourth:
\begin{lstlisting}[language=Mathematica,mathescape,numbers=none,title={Mathematica Session}]
	In[1]:= u=Sqrt[Pi]Cos[Pi I/(4 $\alpha$)] Gamma[-I/(4 $\alpha$)-1/2]/(I/(4 $\alpha$)Gamma[-I/(8$\alpha$)]^2)
	In[2]:= v=-Sqrt[Pi]Cos[Pi I/(4 $\alpha$)] Gamma[I/(4 $\alpha$)-1/2]/(2Gamma[I/(8$\alpha$)+1/2]^2)
	In[3]:= b=2 Sqrt[Pi] Gamma[I/(4 $\alpha$)+3/2](1/Gamma[I/(8$\alpha$)+1/2]^2-1/(Gamma[I/(8$\alpha$)+1]Gamma[I/(8$\alpha$)]))
	In[4]:= bstar=2 Sqrt[Pi] Gamma[-I/(4 $\alpha$)+3/2](1/Gamma[-I/(8$\alpha$)+1/2]^2-1/(Gamma[-I/(8$\alpha$)+1]Gamma[-I/(8$\alpha$)]))
	In[5]:= FullSimplify[u b + v bstar]
	Out[5]= 1
	In[6]:= c=(1/2)/(I/(4 $\alpha$)+3/2)  4/3 Sqrt[Pi] Gamma[I/(4 $\alpha$)+5/2]((1-I/$\alpha$)/Gamma[I/(8$\alpha$)+1/2]^2+(1+I/$\alpha$)/(Gamma[I/(8$\alpha$)+1]Gamma[I/(8$\alpha$)])) (-I $\alpha$/2)
	In[7]:= cstar=(1/2)/(-I/(4 $\alpha$)+3/2) 4/3 Sqrt[Pi] Gamma[-I/(4 $\alpha$)+5/2]((1+I/$\alpha$)/Gamma[-I/(8$\alpha$)+1/2]^2+(1-I/$\alpha$)/(Gamma[-I/(8$\alpha$)+1]Gamma[-I/(8$\alpha$)])) (I $\alpha$/2)
	In[8]:= FullSimplify[Factorial[3](u c + v cstar)]
	Out[8]= -1+I $\alpha$
\end{lstlisting}
This shows us the unique values of $u$ and $v$ that satisfy our initial conditions. 
\end{appendices}

\bibliographystyle{unsrt}
\bibliography{bibliography}

\textit{E-mail address}: \texttt{alexandru.ionut172@gmail.com}

\end{document}